\documentclass[12pt]{article}
\usepackage{epsfig}
\usepackage{graphics}
\usepackage{color}
\usepackage{colordvi}
\usepackage{epsfig}
\usepackage{chicago}

\usepackage{amssymb}

\newcommand{\bea}{\begin{eqnarray}}
\newcommand{\eea}{\end{eqnarray}}

\newcommand{\Pa}{\mbox{pa}}

\newcommand{\dsep}{\, \mbox{\sc ir}_{\delta} }

\newcommand{\cl}{\mbox{\footnotesize cl}}

\newcommand{\aan}{\mbox{\footnotesize An}}

\newcommand{\cale}{{\cal E}}

\newcommand{\lam}{\lambda}
\newcommand{\Lam}{\Lambda}
\newcommand{\bfN}{{\bf N}}

\newcommand{\ld}{\ldots}

\newcommand{\ind}{\bot\!\!\!\bot\,}
\newcommand{\gsep}{\bot\!\!\!\bot_g\,}
\newcommand{\li}{\rightarrow\!\!\!\!\!\!\!/\;\,}
\newcommand{\calf}{{\cal F}}
\newcommand{\calt}{{\cal T}}

\newcommand{\nli}{\rightarrow}
\newcommand{\bsl}{\backslash}

\newcommand{\bo}{\hspace*{\fill}$\Box$}
\newcommand{\nn}{\nonumber}

\newtheorem{defi}{Definition}[section]

\newtheorem{theo}[defi]{Theorem}

\newtheorem{prop}[defi]{Proposition}

\newtheorem{ass}[defi]{Assumption}

\begin{document}

\begin{center}
{\bf \Large Graphical Models for  Marked Point Processes\\[3mm] 
based on Local Independence}
\\
\vspace{8mm}

VANESSA DIDELEZ 

{\em University College London}\\[8mm]
\end{center}

\begin{abstract}
A new class of graphical models capturing the dependence structure of events that
occur in time is proposed. The graphs represent so--called local
independencies, meaning that the intensities of certain types of events are 
independent of some (but not necessarilly all) events in the past.
This dynamic concept of independence is asymmetric, similar to 
Granger non--causality, so that the corresponding local independence graphs 
differ considerably from classical graphical models. Hence a new notion
of graph separation, called $\delta$--separation, is introduced and 
implications for the underlying model as well as for likelihood inference 
are explored. 
Benefits regarding facilitation of
reasoning about and understanding of dynamic dependencies as well as
computational simplifications are discussed.
\end{abstract}


\section{Introduction}

Marked point processes are commonly used to model event history data, 
a term originating from sociology where
it is often of interest to investigate the dynamics behind events such as
finishing college, finding a job, getting married, starting a family, 
durations of unemployment, illness etc.
But comparable data situations also occur in other contexts, e.g.\ 
in survival analyses with intermediate events such as the onset of a side effect,
change of medication etc.\ (cf. Keiding, 1999).
Longitudinal studies and the careful analysis of the underlying processes are 
crucial
for gaining insight into the driving forces of inherently dynamic systems, but
having to deal with the multidimensionality as well as with the dynamic nature
of these systems makes this a very complex undertaking.

Graphical models deal with complex data structures that
arise whenever the interrelationship of 
variables in a multivariate setting is investigated.
Over the last two decades,
they have proven to be a valuable tool for probabilistic modelling and multivariate data 
analysis in such different fields as expert systems and artificial
intelligence (Pearl, 1988; Cowell et al., 1999; Jordan, 1999), 
hierarchical Bayesian modelling (cf.\ BUGS project
{\tt http://www.mrc-bsu.cam.ac.uk/bugs/welcome.shtml}),
causal reasoning (Pearl, 2000; Spirtes et al., 2000), as well as sociological, 
bio--medical, and  econometric applications. For overviews and many 
different applications see for instance the monographs by Whittaker (1990), 
Cox and Wermuth (1996), Edwards (2000). 

While the `classical' graphical models are concerned with representing conditional
independence structures among random variables, variations have been proposed to 
deal with feedback systems (Sprites, 1995; Koster, 1996) but are still based 
on cross--sectional data. The application of graphical models 
to truly time--dependent data, such as event histories or time
series is only slowly getting on its way.
Dahlhaus (2000) has proposed graphical models for multivariate time 
series, but his approach does not capture how the 
present or future of the system depends on or is affected by the past. 
Instead, Eichler (1999, 2000) uses graphs to represent
Granger--causality, which is a dynamic concept of dependence. 
For continuous time one approach, called Dynamic Bayesian Networks, 
is to discretise time and provide directed acyclic graphs that encode the 
independence structure for the transitions from $t$ to $t+1$ (Dean and
Kanazawa, 1989).
The approach proposed by Nodelman et al.\ (2002, 2003) comes closest to
the kind of graphs that we will consider. In their Continuous Time Bayesian
Networks they represent a multistate Markov process with
nodes corresponding to subprocesses and edges corresponding to dependencies
of transition rates on states of other subprocesses.

In this paper, we propose and investigate the properties of graphs that 
represent so--called `local independence' structures in event history data.
The basic idea of local independence is that,
once we know about specific past events, the intensity of a considered 
future event is independent of other past events.
It has been developed by Schweder (1970) for
the case of Markov processes and applied e.g.\ in Aalen et al.\ (1980). 
A generalisation to processes with a Doob--Meyer decomposition can be found in Aalen 
(1987) who focusses on the bivariate case, i.e.\ local dependence
between two processes. 
Here we extend this approach to more than two processes.
The analogy of the bivariate case to Granger-non-causality has been pointed out by
Florens and Foug\`ere (1996), see also Comte and Renault (1996).
Note that the notion of `local independence'
used by Allard et al.\ (2001) is a different one. \\

We first set out the necessary notations and assumptions for marked point processes in 
section 
\ref{sec_mpp} followed by the formal definition of local independence in section
\ref{sec_loc_ind}, the emphasis being on the 
generalisation to a version that allows to condition on the past of other processes
and hence describes dynamic dependencies for multivariate processes.
Section \ref{sec_def} defines graphs that are appropriate to represent local
(in)dependence. 
The main results are given in section \ref{sec_mps}. The properties of local
independence graphs are investigated. In particular we
prove that a new notion of graph separation, called $\delta$--separation,
can inform us about independencies that are preserved after marginalising
over some of the processes. 
In section \ref{sec_lik}, it is shown how the likelihood of a process with given 
local independence graph factorises and implications are discussed.
The potential of local independence graphs is discussed in section \ref{goodfor} 
and proofs are given in the appendix.

\section{Local independence for marked point processes}
Marked point processes are briefly reviewed in section \ref{sec_mpp} 
using the notation of Andersen et al.\ (1993).
In section \ref{sec_loc_ind} the concept of local independence is explained 
in some detail.

\subsection{Marked point processes and counting processes}
\label{sec_mpp}

Let ${\cal E}=\{e_1,\ld, e_K\}$, $K<\infty$, denote the
(finite) \em mark space, \em i.e.\ the set containing all 
\em types \em of events of interest for one observational unit, and
$\cal T$ the \em time space \em in which the observations take place. We assume 
that time is measured continuously so that
we have ${\cal T}=[0,\tau]$ or ${\cal T}=[0,\tau)$ where $\tau<\infty$. 
The \em marked point process (MPP) \em $Y$ consists of \em events \em 
given by pairs of variables $(T_s, E_s)$,
$s=1,2,\ldots,$ on a probability 
space $(\Omega, \calf, P)$ where $T_s\in \cal T$, $0<T_1<T_2
\ld ,$ are the times of occurrences of the respective types of events
$E_s\in {\cal E}$. 
Assume that
the MPP is non--explosive, i.e.\ only a finite number of events occurs in the time
span ${\cal T}$.
The \em mark specific counting processes $N_k (t)$
\em associated with an MPP are then given by
\bea
N_k (t)=\sum_{T_s\leq t}{\bf 1}\{E_s=e_k\},\quad \quad k=1,\ld, K.
\nn
\eea
We write $\bfN=(N_1, \ldots, N_K)$ for the multivariate counting process, and
$\bfN_A$, $A\subset \{1,\ld,K\}$, for the vector $(N_k)_{k\in A}$, 
calling $\bfN_A$ a \em subprocess, \em  with $\bfN_V=\bfN$.

To investigate dependencies of the present on the past it will be important to
have some notation for the history of some subset or all of the processes involved.
Hence we denote the internal filtration of a marked point process by
$\calf_t=\sigma\{(T_s, E_s)\mid T_s\leq t, E_s\in \cale \}$ which is equal to 
$\sigma\{(N_1(s),\ld,N_K(s))\mid s\leq t\}$, while
for $A\subset \{1,\ld,K\}$ we define
the filtrations of a subprocess as $\calf_t^A=\sigma\{\bfN_A(s)\mid s\leq t\}$;
in particular $\calf_t^k$ is the internal filtration of an individual 
counting process $N_k$.

Under rather general assumptions (cf.\ Fleming and Harrington, 1991, p.\ 61), 
a Doob--Meyer decomposition of $N_k(t)$ 
into a compensator and a martingale exists. 
Both these processes depend on the considered filtration
which is here taken to be
the internal filtration of the \em whole \em MPP $Y$. We will assume throughout 
that
all the $\calf_t$--compensators $\Lam_k$ are absolutely continuous and predictable
so that 
intensity processes $\lam_k(t)$ exist, which are taken to be predictable versions 
of the derivatives of the compensators, i.e.\
$\Lam_k(t)=\int_0^t \lam_k(s)ds$. Heuristically we have (Andersen et al., 1993, 
p.\ 52)
\bea
\lam_k(t)dt=E(N_k(dt)\mid \calf_{t^-}) \label{compens}.
\eea
More formally this means that the 
differences $N_k-\Lam_k$ are $\calf_t$--martingales.

An interpretation of property (\ref{compens}) is
that given the information on the history of the whole MPP up to just
before time $t$, $\lam_k(t)dt$ is our best prediction of the immediately following
behaviour of $N_k$. Note that for the
theory developed in this paper the setting can slightly be generalised, 
not requiring absolute continuity of compensators (cf.\ Didelez, 2000), as might
be relevant when certain types of events can only occur at fixed times.

It will be important to distinguish between
the $\calf_t$--intensity based on the past of the \em whole \em 
MPP, and the
$\calf_t^A$--intensities based on the past
of the subprocess on marks in $A$. The latter
can be computed using the innovation theorem (Br\'emaud, 1981, pp.\ 83),
and a way of doing so, especially relevant to our setting, 
is given in Arjas et al.\ (1992).

The following is a standard assumption in counting process theory but we want to
highlight it as it plays a particularly important role for local independence
graphs.

\begin{ass} No jumps at the same time \em \label{ass_jumps} \\
The  ${\cal F}_t$--martingales 
$N_k-\int \lam_k(s) ds$ are assumed to be orthogonal for $k\in\{1,\ld,K\}$,
meaning that none of $N_1, \ldots, N_K$ jump at the
same time. This is implied
by the above assumption that all compensators are absolutely continuous if in 
addition no
two counting processes $N_j$ and $N_k$ are counting the same type of event. 
\end{ass}

Assumption \ref{ass_jumps}
might be violated, e.g.\ when investigating the survival times of
couples and there is a small but non--zero chance that they die at the same
time, in a car accident for instance. The reason for imposing this assumption
is that we want to explain dependencies between events by the past not by
common innovations. If one wants to allow events to occur at the same time, 
then such
a simultaneous occurrence defines a new mark in the mark space.

Note that general multi--state processes can be represented as
marked point process with every transition between two states being a mark.
This is explored in more detail in Didelez (2007) for Markov processes.

\subsection{Local independence}
\label{sec_loc_ind}

The bivariate case is defined as follows (cf.\ Aalen, 1987).

\begin{defi} Local independence (bivariate). \label{def_li_biv} \\
\em
Let $Y$ be an MPP with $\cale=\{e_1, e_2\}$ and $N_1$ and $N_2$ the associated 
counting processes on $(\Omega, \calf, P)$. 
Then, $N_1$ is said to be \em locally independent \em of $N_2$ over $\cal T$ if
$\lam_1(t)$ is measurable w.r.t\ ${\cal F}^1_t$ for all $t\in \calt$. 
Otherwise we speak of local dependence. 

The process $N_1$ being locally independent of $N_2$ is symbolised by $N_2\li N_1$. 
Interchangeably we will sometimes say that 
$e_1$ is locally independent of $e_2$, or $e_2\li e_1$.
\end{defi}

The essence of the above definition is that the intensity  $\lam_1(t)$, i.e.\ our
`short--term' prediction of $N_1$, remains the same under the reduced filtration $\calf_t^1$ as
compared to the full one $\calf_t$. This implies that we do not lose any essential information
by ignoring how often and when event $e_2$ has occurred before $t$.
One could  say that if $N_2\li N_1$ then the presence of $N_1$ is conditionally
independent of the past of $N_2$ given the past of
$N_1$, or heuristically\footnote{Note that (\ref{condind1}) is an informal way of 
saying that $N_1(dt)$ is conditionally independent of $\{T_s|T_s<t,E_s=2, s=1,2,\ldots\}$
given $\{T_s|T_s<t,E_s=1, s=1,2,\ldots\}$. This and similar statements later, like 
(\ref{condind2}), (\ref{part_fact}) or (\ref{ind_1}), 
should be interpreted correspondingly.}
\begin{eqnarray}
N_1(t)\ind {\cal F}_{t^-}^{2}|{\cal F}_{t^-}^1,
\label{condind1}
\end{eqnarray}
where $A\ind B|C$ means `$A$ is conditionally independent of $B$ given $C$'
(cf.\ Dawid, 1979). For general processes, this is
a stronger property than local independence but it holds for marked point processes 
with Assumption \ref{ass_jumps} as 
their distributions are determined by the intensities. 
Note that (\ref{condind1}) does not imply that for $u>0$:
$N_1(t+u)\ind {\cal F}_{t^-}^{2}|{\cal F}_{t^-}^1$ --- hence the name
\em local \em independence. Also,
$N_1(t)\ind N_2(t)$ will only hold if the two processes
are mutually locally independent of each other.
Without Assumption \ref{ass_jumps} the $\calf_t^1$--measurability of 
$\lam_1(t)$ would, for instance, trivially be true if
$e_1=e_2$ but, in such a case, we would not want to speak of 
independence of $e_1$ and $e_2$. \\

\noindent {\bf Example 1.} \em Skin disease. \em \\
In a study with women of a certain age,
Aalen et al.\ (1980) model two events in the life of an individual women: 
occurrence of a 
particular skin disease and onset of menopause. Their analysis reveals that
the intensity for developing this skin disease is greater once the menopause 
has started than before. In contrast, and as one would expect, the intensity for 
onset of menopause
does not depend on whether the person has earlier
developed this skin disease. We can therefore
say that menopause is locally independent of this skin disease but not
vice versa. Note that, in whatever way the onset of skin disease and menopause 
are measured, it is assumed that they do not
start systematically at exactly the same time, corresponding to the above `no 
jumps at the same time' assumption.\\

Let us now turn to the case of more than two types of events. 
This requires conditioning on the past of other processes as follows.

\begin{defi} Local independence (multivariate). \label{def_li_mul}\\
\em
Let $\bfN=(N_1, \ldots, N_K)$ be a multivariate counting processes associated with an MPP. Let
further $A,B,C$ be disjoint subsets of $\{1,\ld,K\}$. We then say that a \em subprocess
$\bfN_B$ is locally independent of
$\bfN_A$ given $\bfN_{C}$ \em over $\cal T$
if all $\calf_t^{A\cup B\cup C}$--intensities $\lam_k$, $k\in B$, are measurable
with respect
to $\calf_t^{B\cup C}$ for all $t\in {\cal T}$. This is denoted by $\bfN_A\li \bfN_B \mid \bfN_{C}$
or briefly $A\li B \mid C$.
Otherwise, $\bfN_B$ is \em locally dependent \em on $\bfN_A$ given $\bfN_C$,
i.e.\ $A\nli B \mid C$. If $C=\emptyset$ then $B$ is \em marginally \em locally
(in)dependent of $A$.
\end{defi}

Conditioning on a subset $C$ thus means that we retain the information about
whether and when events with marks in $C$ have occurred in the past when
considering the intensities for marks in $B$. If these intensities are 
independent of the information on whether and when events with marks in $A$ 
have occurred
we have conditional local independence.
In analogy to (\ref{condind1}) the above definition of multivariate local independence
is with Assumption \ref{ass_jumps} equivalent to
\begin{eqnarray}
N_B(t)\ind {\cal F}^A_{t^-}|{\cal F}_{t^-}^{B\cup C} \quad \forall\, t\in \calt.
\label{condind2}
\end{eqnarray}

\noindent {\bf Example 2.} \em Home visits \em \\
This example is not taken from the literature but is inspired by real studies 
(e.g.\ Vass et al., 2002, 2004). In some countries programmes exist to assist 
the elderly  through regular home visits by a nurse. This is meant to reduce 
unnecessary hospitalisations while increasing the quality of life for the person. 
It is hoped that such a programme
increases the survival time.  The
times of the visits as well as the times and durations of hospitalisation are
monitored. In addition, it is plausible that 
the underlying health status of the elderly person 
may also affect the rate of hospitalisation and predict 
survival. This interplay of events for an individual elderly person
can be represented as an MPP if `health status' is regarded as a
multistate process. Assume that the timing of the
home visits is determined externally, e.g.\
by the availability of nurses which has nothing to do with the patient's 
development, i.e.\ the visits are assumed locally independent of
all the remaining processes. 
It might then be of interest to investigate whether the visits affect only the 
rate of hospitalisation directly, i.e.\ whether survival is locally independent
of the visits process given the hospitalisation and health history or even given
only a subset thereof.\\

As can easily be checked, local (in)dependence needs to be neither symmetric, 
reflexive nor transitive. However, since in most practical situations 
a subprocess depends at least on its own past we will assume throughout
that local dependence is reflexive. An example for a subprocess that depends
only on the history of a different subprocess and not on its own history
is given in Cox and Isham (1980, p.\ 122).\\

In order to see the relation with local independence, we 
briefly review Granger non--causality
(Granger, 1969). Let $X_V=\{X_V(t)\mid t\in {\bf Z}\}$ with $X_V(t)=(X_1(t),
\ldots, X_K(t))$ be a multivariate time series, where $V=\{1,\ldots, K\}$ is 
the index set. For any $A\subset V$ we define $X_A=\{X_A(t)\}$ as the multivariate
subprocess with components $X_a$, $a\in A$. Further let $\overline X_A(t)=
\{X_A(s)\mid s\leq t\}$. Then, for disjoint subsets $A,B\subset V$ we say that
$X_A$ is strongly Granger-non causal for $X_B$  if 
\[
X_B(t)\ind \overline X_A(t-1) \mid \overline X_{V\backslash A}(t-1),
\]
for all $t \in Z$. The interpretation is similar as
for local independence, i.e.\ the present value of  $X_B$ is independent of the
past of $X_A$ given its own past and the one of all other components 
$C=V\backslash (A \cup B)$, in
analogy to (\ref{condind2}). Also note that the above does not imply that
$X_B(t+u)\ind \overline X_A(t-1) \mid \overline X_{V\backslash A}(t-1)$ for $u>0$, 
again
analogous to local independence. Eichler (1999, 2000) investigates a graphical 
representation and rules to determine when the condition $\overline X_{V\backslash 
A}(t-1)$ can be reduced to proper subsets $\overline X_{C}(t-1)$, 
$C\subset V\backslash A$.\\

Finally, let us indicate how the definition of local
independence can be generalised to stopped processes. This is relevant
when there are absorbing states such as death. In that case
all other events will be locally dependent on this one because all intensities
are zero once death has occurred. However, the dependence is `trivial' and not
of much interest. Let $T$ be a ${\cal F}_t$--stopping time
and let $\bfN^T=(N^T_1, \ldots, N^T_K)$ be the
multivariate counting process stopped at $T$. Then the intensities of $N^T_k$
are given by $\lam_k^T$, $k\in V$, and local independence can be formalised as
follows.

\begin{defi} Local independence for stopped processes. \label{def_stopped} \\
\em
Let $\bfN^T=(N_1^T, \ldots, N_K^T)$ be a multivariate counting processes 
associated to an MPP and stopped at time $T$. 
Then we say that $A\li B \mid C$ 
if there exist 
$\calf_t^{B\cup C}$--measurable  processes $\tilde \lam_k$, $k\in B$,
such that the 
$\calf_t^{A\cup B\cup C}$---intensities of $\bfN^T_B$ are given by
$\lam_k^T(t)=\tilde \lam_k(t) {\bf 1}\{t\leq T\}$, $k\in B$.
\end{defi}

The local independencies in a stopped process have to be interpreted as
being valid as long as $t\leq T$.

\section{Local independence graphs}\label{sec_li_graphs}

We first give the definition of local independence graphs
and then investigate what can be read off these graphs.

\subsection{Definition of local independence graphs}\label{sec_def}

An obvious way of representing the local independence structure of an MPP by a graph is to
depict the marks as vertices and to use an arrow as symbol for local dependence as in
the following small example. \\

\begin{figure}[h]
\begin{center}
\leavevmode
\epsfclipon
\epsfxsize55mm
\epsffile{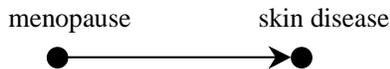}
\epsfclipoff
\caption{Local independence graph for skin disease example.}\label{fig_1}
\end{center}
\end{figure}

\noindent {\bf Example 1 ctd.} \em Skin disease. \em\\
The local independence graph for the relation between menopause and skin disease
is very simple, cf.\ Figure \ref{fig_1}. 
Notice that even with this simple example there is no way of expressing the 
local independence using a classical graph based on conditional
independence for the two times $T_1=$ `time of occurrence of skin disease' and
$T_2=$ `time of occurrence of menopause' as these are simply dependent.\\

For general local independence structures, we will have directed graphs that 
may have more than one directed edge between a pair of vertices, in case
of mutual local dependence, and that may have cycles. More formally,
a \em graph \em is a pair $G=(V,E)$, where $V=\{1,\ld, K\}$  is a finite set of vertices and
$E$ is a set  edges.
The graph is said to be 
\em directed \em if $E\subset  \{(j,k)|j,k\in V, j\not=k\}$.
Later we will also need the notion of an \em undirected \em graph where 
$E\subset  \{\{j,k\}|j,k\in V, j \not=k\}$.
Undirected edges $\{j,k\}$ are depicted by lines, $j$ ---
$k$, and directed edges $(j,k)$ by arrows, $j\longrightarrow k.$ If $(j,k)\in E$ and $(k,j)\in
E$ this is shown by stacked arrows, $j \stackrel{ \mbox{ \normalsize $\longrightarrow$ }}{
\longleftarrow} k$. \\

The following property (\ref{pair}) is called the \em pairwise dynamic 
Markov property\em , where we say \em dynamic \em  
to emphasise the difference to graphs based on conditional independence.

\begin{defi} Local independence graph\label{defi_pair}\\
\em
Let $\bfN_V=(N_1,\ldots,N_K)$ be a multivariate counting process associated to an MPP $Y$ with
mark space $\cale=\{e_1, \ldots , e_K\}$. Let further $G=(V,E)$ be a directed graph,
$V=\{1,\ld,K\}$. Then, $G$ is called a \em local independence graph \em of $Y$  if
\bea
\mbox{for all }j,k\in V: \quad
(j,k)\notin E \quad
\Rightarrow \quad \{j\}\li \{k\}|V\bsl \{j,k\} \label{pair}.
\eea
\end{defi}
\vspace{5mm}

\begin{figure}[h]
\begin{center}
\leavevmode
\epsfclipon
\epsfxsize125mm
\epsffile{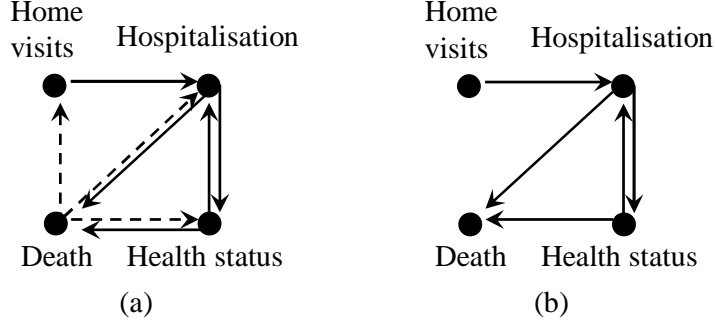}
\epsfclipoff
\caption{Home visits example; (a) for the whole process;  
(b) for the stopped process.}
\label{dementia3}
\end{center}
\end{figure}

\noindent {\bf Example 2 ctd.} \em Home visits \em \\
The graph in Figure \ref{dementia3}(a) is for the whole process, while 
(b) shows the local independencies for the stopped process (stopping when 
death occurs). There are no arrows into `Home visits' representing that the
rate of visits is locally independent of `Hospitalisation' given `Health status'
as well as of `Health status' given `Hospitalisation' (while the person is
still alive), reflecting that the visits are determined externally. The latter local 
independence might be violated if the nurses, on their own account,
increase the frequency of their visits when they notice that the person's health is
deteriorating. The graph further represents that survival is locally 
independent of the visits given hospitalisation and health history, and that
the health process is also locally independent of the visits given hospitalisation 
history (while the person is still alive obviously). These absent edges could
reflect the null--hypothesis when investigating whether the visits affect 
survival in other ways than through changing the rate of hospitalisation.

\subsection{Dynamic Markov properties}\label{sec_mps}

The local independence graphs as defined above allow (under mild assumptions) 
more properties to be
read off, concerning the dependence structure, than just those given
by (\ref{pair}). We may in particular query the graph with the aim of dimension
reduction, i.e.\ with questions about which other processes can be ignored while 
investigating certain local independencies. The \em local dynamic Markov property
\em addressed in section \ref{sec_local} tells us about the immediately 
relevant information when considering a single mark $e_k$ and corresponding $N_k$.
Further, the \em global dynamic 
Markov property \em in section \ref{sec_global}  gives graphical rules to identify 
when the \em separating \em set
itself can be reduced, i.e.\ when in (\ref{pair}) we do not need to condition on
\em all \em $V\bsl \{j,k\}$ but just on a true subset. For this we need the notion
of $\delta$--separation also introduced in section \ref{sec_global}.\\

Some more graph notation will be required.
A \em path \em between two nodes is defined in the obvious way (the formal definition is given in
Appendix \ref{app_delta}): we distinguish between \em undirected paths \em for undirected graphs,
\em directed paths, \em preserving the direction of edges, for directed graphs 
and \em trails \em
for connections in directed graphs that do not preserve the direction.
For directed graphs we further require the following almost self-explanatory 
notations. If $a\longrightarrow b$ then $a$ is called a \em parent \em of $b$ and
$b$ is a \em child \em of $a$ (if 
$a \stackrel{ \mbox{ \normalsize $\longrightarrow$ }}{
\longleftarrow} b$ then $a$ is both, a child and a parent of $b$); 
pa$(A)$ denotes the set of all parents of nodes
in $A\subset V$ without $A$ itself, and  ch$(A)$ analogously the set of children of $A$.
The set cl$(A)=\Pa(A)\cup A$ is called the \em closure \em of $A$.
If there is a directed path from $a$ to $b$ then $a$ is an \em ancestor \em of
$b$ and $b$ is a descendant of $a$; the corresponding set notation is an$(A)$ and 
de$(A)$ (always excluding $A$ itself). Consequently, nd$(A)=V\backslash
($de$(A)\cup A)$ are the \em non--descendants \em of $A$.
If pa$(A)=\emptyset$, then $A$ is called \em ancestral. \em In general, An$(A)$
is the \em smallest ancestral set \em containing $A$, given by $A$ $\cup$ an$(A)$.

\subsubsection{Local dynamic Markov property}\label{sec_local}

\begin{defi} Local dynamic Markov property\\
\em
Let $G=(V,E)$ be a directed graph. For an MPP $Y$ the property
\bea
\mbox{for all } k\in V:\quad V\bsl \mbox{cl}(k) \li \{k\}\mid \mbox{pa}(k), \label{locmp}
\eea
is called the \em local dynamic Markov property w.r.t.\ $G$.
\end{defi}

In other words, property (\ref{locmp}) says that every $\calf_t$--intensity 
$\lam_k$ is $\calf_t^{\cl(k)}$--measurable, which clearly implies
that for any ancestral set $A$ the intensity $\lam_A$ is $\calf_t^{A}$--measurable. 
This property could for instance be violated if two components in pa$(k)$
were a.s.\ identical 
which is however  prevented by the orthogonality assumption \ref{ass_jumps}.
As shown in  Appendix \ref{app_proof},
the exact condition for the property (\ref{locmp}) to follow from 
(\ref{pair}) is that
\bea
{\cal F}_t^{A}\cap{\cal F}_t^{B}={\cal F}_t^{A\cap B}
\quad \forall \; A,B\subset V,
\quad \forall\, t\in {\cal T} \label{filter},
\eea
where we define ${\cal F}^{\emptyset}=\{\emptyset, \Omega\}$.
Property (\ref{filter}) is called `conditional measurable separability' (Florens et 
al., 1990) and formalises the intuitive notion that the components of $\bfN$ are
`different' enough to ensure that common events are necessarily due to 
common components. \\

\noindent {\bf Example 2 ctd.} \em Home visits \em \\
Let us consider the question whether the four processes are sufficiently
different to ensure (\ref{filter}). If the health process is
measured in a way such that it is determined by the number and duration of past
hospitalisation, not taking any other information into account, this assumption
might be violated. However, it makes sense and we will assume for this example that
the `Health status' reflects more aspects of a person's health than just past
hospitalisations. Then it seems plausible that (\ref{filter}) is satisfied as the
other processes are clearly capturing different information anyway. 
Consequently we can
use the local dynamic Markov property to read off that the visits process
is locally independent of both, hospitalisation and health status (while the 
person is still alive).

\subsubsection{$\delta$--separation and the global dynamic Markov property}\label{sec_global}

In undirected graphs we say that subsets $A,B\subset V$ are separated by $C\subset V$ if any
path between elements in $A$ and elements in $B$ is intersected by $C$. This is symbolised by
$A\gsep B|C$. In classical graphical models every such separation induces 
conditional independence between $A$ and $B$ given $C$ 
regardless of whether $(A,B,C)$ is a partition
of $V$ or not. This can obviously lead to considerable dimension reduction if
$C$ is chosen minimally and the graph is sparse. To obtain a similar result for
local independence graphs we require a suitable notion of separation called
$\delta$--separation, introduced below after some more graph notation.\\

The \em moral \em  graph $G^m$ is given by
inserting undirected edges between any two vertices that have a common child (if they are not
already joined) and then making all edges undirected (two directed edges 
between a pair of nodes are replaced by one undirected edge). This procedure of 
moralisation will  also be applied to an \em induced subgraph \em $G_A$, $A\subset V$, 
defined as $(A,E_A)$ with $E_A$ the subset of $E$ containing
only edges between pairs of nodes in $A$. Finally, for $B\subset V$, let $G^B$ denote 
the graph obtained by deleting all directed edges of
$G$ starting in $B$.

\begin{defi} $\delta$--Separation. \label{defi_delta_sep}\\
\em
Let $G=(V,E)$ be a directed  graph. 
Then, we say for pairwise disjoint subsets $A,B,C\subset V$ that $C$ \em $\delta$--separates
$A$ from $B$ in $G$ \em if
$A\gsep B|C$ in the undirected graph $(G^B_{\aan(A\cup B\cup C)})^m$
(the case of non--disjoint $A,B,C$ is given in Appendix \ref{app_delta}).
\end{defi}

Note that except for the fact that we delete edges starting in $B$,  which makes
$\delta$--separation asymmetric, the definition parallels the one for DAGs. 
This initial edge deletion can heuristically be explained by the fact that we
want to separate the present of $B$ from the past of $A$ and hence we disregard
the `future' of $B$ which is where the edges out of $B$ point to; for the same reason only the
ancestral set An$(A\cup B\cup C)$ is considered. As for DAGs the
insertion of moral edges is necessary whenever we condition on a common `child'
due to a `selection effect' by which two marginally independent variables (or processes) that
affect a third variable (or process) become dependent when conditioning on this third variable.
Further properties of $\delta$--separation are discussed in
Didelez (2006).

\begin{defi} Global dynamic Markov property\\
\em
Let $\bfN_V=(N_1,\ldots,N_K)$ be a multivariate counting process associated to an MPP $Y$ and
$G=(V,E)$ a directed graph. The property that
\bea
\mbox{for all disjoint }
A,B,C\subset V: \quad
C\; \delta\mbox{--separates }A\mbox{ from }B\mbox{ in }G\quad
\Rightarrow \quad A\li B\mid C. \label{glob}
\eea
is called the \em global dynamic Markov property w.r.t.\ $G$.
\end{defi}

The significance of the global Markov property is that it provides a way
to verify whether a subset $C\subset V\bsl (A\cup B)$ is given such that
$A\li B|C$, i.e.\ local independence is preserved even \em when ignoring \em 
information on the past of processes in $V\bsl(A\cup B\cup C)$. Of course, (\ref{glob})
is only meaningful if it can be linked to the definition of local independence graphs
addressed next.

\begin{theo} Equivalence of dynamic Markov properties \label{theo_equiv}\\
\em
Let $Y$ be a marked point process and $G=(V,E)$ a directed graph. 
Under the assumption
of (\ref{filter}) and further regularity conditions (cf.\ Appendix \ref{app_proof}), 
the pairwise, local and
global dynamic Markov properties, i.e.\ 
(\ref{pair}), (\ref{locmp}) and (\ref{glob}) are equivalent.
\end{theo}
The proof is given in Appendix \ref{app_proof}. \\

\begin{figure}[h]
\begin{center}
\leavevmode
\epsfclipon
\epsfxsize125mm
\epsffile{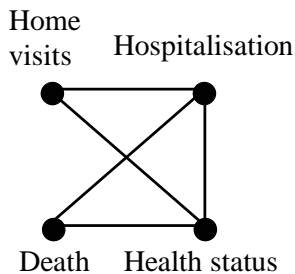}
\epsfclipoff
\caption{Moral graph for home visits example.}
\label{dementia2}
\end{center}
\end{figure}

\noindent {\bf Example 2 ctd.} \em Home visits. \em\\
The underlying health status of an elderly person may be difficult to measure
accurately in practice. Let us therefore investigate the local independence
structure when ignoring this underlying process altogether, in particular
consider the question of whether from Figure \ref{dementia3} we can infer
that survival is still locally independent of `Home visits' 
given only the hospitalisation but 
ignoring the health process. Graphically 
this means we have to check whether the node `Hospitalisation' alone
separates `Home visits' from `Death'. As can be seen from the corresponding 
moral graph (for the stopped process) in Figure \ref{dementia2} this is not the case.
Hence, even though the home visits are assumed to be determined externally in
Figure \ref{dementia3} and do not affect survival directly, ignoring the 
underlying health process may lead to a `spurious' local dependence of survival
on the home visits. The reason is that, for instance,
a history of hospitalisation with preceding
home visit predicts  survival differently from a hospitalisation without preceding
home visit --- the former might mean that the health was especially bad and hence
hospitalisation was necessary, while the latter allows minor health problems that
could have been treated by a nurse who was not available. \\[4mm]

Notice that intuitively it is clear that if the intensity of the visits
depended on the underlying health status, i.e.\ if there was a directed edge from
`Health' to `Visits', we could talk of confounding. Hence it is rather surprising
that even when the frequency of the visits is controlled externally we may find
a spurious dependence. For the discrete time case Robins (1986, 1997) has 
demonstrated that nevertheless in a situation like Figure \ref{dementia3} 
we can draw causal conclusions
even when no information on the underlying health process is available. 
However, standard methods that just model the intensity for survival 
with time varying  covariates for 
the times of previous home visits and  hospitalisations will
typically give misleading results due to the conditional association
between `Home visits' and `Death' given `Hospitalisation'.\\

\begin{figure}[h]
\begin{center}
\leavevmode
\epsfclipon
\epsfxsize125mm
\epsffile{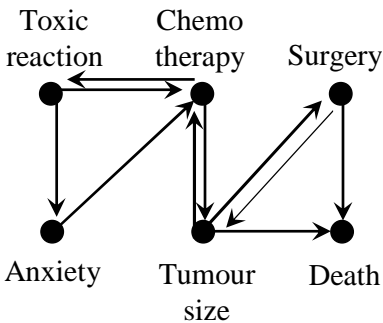}
\epsfclipoff
\caption{Local independence graph for chemo therapy example (stopped process).}
\label{chemo1}
\end{center}
\end{figure}

\noindent {\bf Example 3.} \em Chemo therapy cycles \em \\
To see a more complex example of a local independence graph consider a 
hypothetical study (which is inspired by real studies) 
where early stage breast cancer 
patients are observed over a period of several months during which they 
receive at least one but usually more
cycles of chemo therapy. The tumour size is monitored through palpation.
The doctors will consider removing the tumour by surgery if the size does not
decrease and surgery is almost certain if it increases. 
Furthermore, the chemo therapy may be 
delayed or discontinued if the patient shows a toxic reaction following the
treatment but also if the patient requests a delay, possibly due to an increased
state of anxiety. Except for tumour size all processes count one type of event that
can occur once or more often. Tumour size is measured categorically depending on the number 
and palpable size of lesions and can be regarded as a multistate process.
Figure \ref{chemo1} shows a hypothetical  local independence structure. 
For instance it assumes that survival
locally depends on the tumour size and whether surgery has taken place, but once
this information is given none of the other processes are relevant for the
intensity of death. Note that this particular assumption could plausibly
be violated because toxic reactions and anxiety may reflect other health
problems, but for simplicity we will assume that all patients are `healthy' 
except for the breast cancer so that this violation is excluded.

\begin{figure}[h]
\begin{center}
\leavevmode
\epsfclipon
\epsfxsize125mm
\epsffile{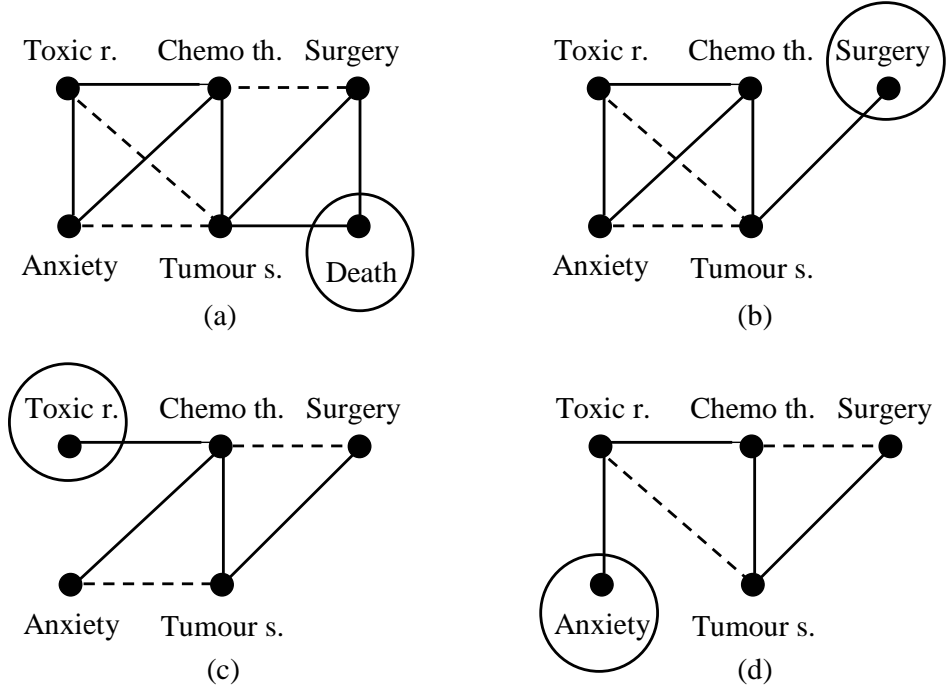}
\epsfclipoff
\caption{Different moral graphs for chemo therapy example. Dotted edges have been
added due to a `common child' in Figure \ref{chemo1} and circles indicate
that arrows out of these nodes have been deleted before moralising.}
\label{chemo2}
\end{center}
\end{figure}

Figure \ref{chemo2} shows the different moral graphs constructed from Figure
\ref{chemo1} to investigate $\delta$--separations. Graph (a) shows 
$\delta$--separation from the node `Death' allowing us to read off, for instance, 
that `Chemo therapy' is not $\delta$--separated from `Death' by `Tumour size' alone
reflecting that chemo therapy predicts survival if surgery history is ignored. This 
is plausible because knowing that for example a decrease in tumor size
was preceded by a treatment cycle is informative for surgery, making it less 
likely than without preceding chemo therapy; and whether surgery has taken place, 
in turn, predicts the survival chances. As `Anxiety' is
problematic to observe and measure
we may further be interested in the question of when it 
can be ignored. We see from the graph (a) that `Death' is locally independent of
`Anxiety' given either the set
$\{$`Surgery', `Tumour size'$\}$ or $\{$`Chemo therapy', `Tumour size'$\}$, the 
latter implying that once we know the chemo therapy history in addition to the
development of the tumour size then anxiety will not inform
us any further about the intensity for death regardless of whether surgery 
and toxic reaction history is known or not. But note that even though `Anxiety'
does not affect `Tumour size' directly the latter has to be part of the separating
set. Figure \ref{chemo2}(b) shows that `Anxiety' is $\delta$--separated from `Surgery'
by any set that includes `Tumour size' and similarly for (c) that it is $\delta$--separated
from `Toxic reaction' by any set that includes `Chemo therapy' --- in these two
cases $\delta$--separation does not tell us more than the local dynamic Markov
property (\ref{locmp}). From graph (d) we see that `Anxiety' itself is locally independent 
of `Chemo therapy'
and `Tumour size' given `Toxic reaction' and of `Surgery' given either 
`Toxic reaction' or the set $\{$`Chemo therapy', `Tumour size'$\}$.

\subsection{Likelihood factorisation and implications}\label{sec_lik}

In order to discuss properties and implications for the likelihood for graphical 
MPPs we will regard the data consisting of times and types of events
$(t_1,e_1), (t_2, e_2),\ldots, (t_n,e_n)$ as a realisation of the 
\em history process \em $H_t=\{(T_s, E_s)|T_s\leq t\}$. As for
filtrations, $H_{t^-}$ denotes the \em strict pre--$t$ history process. \em 
Additionally,
$H_t^A$, $A\subset \{1,\ld,K\}$, defined as
\[
H_t^A=\{(T_s, E_s)\mid T_s\leq t \mbox{ and } \exists\, k\in A : E_s=e_k, s=1,2\ld \}
\]
denotes \em the history process restricted to the marks in $A$. \em Any set of marked points
for which it holds that $t_s=t_u, s\not= u,$ implies $e_s=e_u$ can be a history, i.e.\ a
realization of $H_t$. Note that (up to completion by null sets) 
the different filtrations can be regarded as being generated by
the history processes, i.e.\ $\calf_t^A=\sigma\{H_t^A\}$, $A\subset \{1,\ld,K\}$.

Before deriving the likelihood for a given local independence graph, 
we recall it  for the general case. Based on the mark specific 
intensity processes $\lam_k(t)$ the corresponding \em crude \em 
intensity process is given by $\lam(t)=\sum_{k=1}^K\lam_k(t)$.
This is the intensity process of the cumulative counting process 
$\sum_{k}N_k$.
The likelihood process $L(t| H_t)$ is then
given as 
\bea
L(t| H_t)=\prod_{T_s\leq t}\lam_{E_s}(T_s)\, \cdot\, \exp \left( -\int_0^t\lam(s)\,ds \right).
\label{cont_lik}
\eea
To see how the likelihood is affected
by $G$ being a local independence graph of $Y$, 
we first rewrite (\ref{cont_lik}) as follows:
\bea
L(t| H_{t})
&=& \prod_{k=1}^K \prod_{T_s\leq t}\lam_k(T_s)^{\mbox{\footnotesize \bf 1}\{E_s=e_k\}} \,
\cdot\,
\exp \left( -\int_0^t\sum_{k=1}^K \lam_k(s)\,ds \right) \nn \\[2mm]
&=& \prod_{k=1}^K \left(
\prod_{T_{s(k)}\leq t}\lam_k(T_{s(k)}) \, \cdot\, \exp \left( -\int_0^t\lam_k(s)\,ds \right)
\right),
\nn
\eea
where $T_{s(k)}$ with $E_{s(k)}=e_k$ are the occurrence times of mark $e_k$. 
The inner
product of the above can be regarded as the mark specific likelihood and
is denoted by $L_k(t| H_{t})$. Now, 
by the definition of a local independence graph and
the equivalence of the pairwise and local dynamic Markov properties under condition
(\ref{filter}) we have that $\lam_k(s)$ is $\calf_t^{\cl(k)}$--measurable,
where cl$(k)$ is the closure of node $k$. Hence, it follows that
\bea
L_k(t| H_{t})
&=& L_k(t| H_{t}^{\cl(k)}),
\eea
i.e.\ the mark specific likelihood 
$L_k$ based on the whole past remains the same if the available
information is restricted to how often and when those marks that are parents of $e_k$ in the graph and
$e_k$ itself have occurred in the past, symbolised by $H_{t}^{\cl(k)}$.

It follows that under (\ref{filter}) the likelihood  factorises as
\bea
L(t| H_{t}) = \prod_{k\in V} L_k(t| H_{t}^{\cl(k)}), \label{likfact}
\eea
which parallels the factorisation for DAGs where the joint density is decomposed into the
univariate conditional distributions given the parents. Here, we replace the parents by the
closure because we also need to condition on the past of a component itself which is not
required in the static case.\\

\noindent {\bf Example 2 ctd.} \em Home visits \em\\
From Figure \ref{dementia3}(b) we obtain the following factorisation
\[
L(t| H_{t}) = L_{vi}(t| H_{t}^{\{vi\}})L_{ho}(t| H_{t}^{\{vi, ho, hs\}})
L_{hs}(t| H_{t}^{\{ho, hs\}})L_{d}(t| H_{t}^{\{d,ho,hs\}})
\]
for $t\leq$ time of death, where $vi$, $ho$, $hs$, $d$ stand for `visits',
`hospitalisation', `health status' and `death' respectively.\\

Two consequences of the above factorisation regarding the relation of local and
conditional independence are given next.

\begin{theo} \label{theo_cis} Conditional independencies \\ \em
For an MPP with local independence graph $G$ and disjoint
$A,B,C\subset V$, such that $C$ separates $A$ and $B$, i.e.\ $A \gsep B \mid C $, 
in $(G_{\aan(A\cup B\cup C)})^m$, we have
\bea
\calf_t^{A}\ind \calf_t^{B}\mid \calf_t^{C} \quad \forall\,t \in \calt. 
\label{total_fact}
\eea
\end{theo}

The proof is given in appendix \ref{app_cis}.
The graph separation  that $A,B,C$ have to satisfy for (\ref{total_fact}) 
implies
that for each $k\in C$ the $\calf_t^{A\cup B\cup C}$--intensity $\lam_k$
is either $\calf_t^{A\cup C}$-- or $\calf_t^{B\cup C}$--measurable, otherwise
$C$ could not separate $A$ and $B$ in the moral graph.
Also, of course, we have that for each 
$k\in A$ (resp.\ $B$) the $\calf_t^{A\cup B\cup C}$--intensity
is $\calf_t^{A\cup C}$
(resp. $\calf_t^{B\cup C}$)--measurable; otherwise there would be 
edges linking $A$ and $B$ in the graph and they could not be separated.
A property similar to (\ref{total_fact}) has been noted by Schweder (1970, Theorems 3 and 4) for Markov
processes. With a similar argument we can reformulate (\ref{condind1}) and (\ref{condind2}): 
for any $B\subset V$ we have
\bea
\bfN_B(t)\ind \calf_{t^-}^{V\bsl \cl(B)} \mid \calf_{t^-}^{\cl(B)},
\label{part_fact}
\eea
i.e.\ the present of $\bfN_B$ is independent of the past of
$\bfN_{V\bsl \cl(B)}$ given the past of $\bfN_{\cl(B)}$.\\

\noindent {\bf Example 3 ctd.} \em Chemo therapy cycles \em \\
Let $A=$`Surgery', $B=$`Toxic reaction' and
$C=\{$`Chemo therapy', `Tumour size'$\}$. Then $(G_{\aan(A\cup B\cup C)})^m$
is the same as Figure \ref{chemo2}(a), where the node
`Death' could be omitted as we are conditioning on the patient being alive anyway,
and indeed $A$ and $B$ are separated by $C$. With 
(\ref{total_fact}) we can infer that at any time $t$ 
(before death) the whole surgery history, i.e.\ whether and when surgery 
has taken place before $t$, is independent of whether and when toxic reactions
have occured given we know the tumour size development up to $t$
and when chemo therapy has been administered. 
Note that as mentioned above we have
for the nodes in $C$ that the $\calf_t^{A\cup B\cup C}$--intensity for 
`Tumour size' is  $\calf_t^{A\cup C}$-- and  the one for `Chemo therapy' is 
$\calf_t^{B\cup C}$--measurable; the latter can be seen, using 
$\delta$--separation, by checking that
`Chemo therapy' is locally independent from $A=$`Surgery' given $\{$`Chemo
therapy', `Toxic reaction' $\}$  (the relevant moral
graph happens to be the same as Figure \ref{chemo2}(b)).

\subsection{Extensions}\label{sec_ext}

Local independence graphs can easily be extended to include time--fixed
covariates, such as sex, age and socioeconomic background of patients etc.
The filtration at the start, ${\cal F}_0$, then has to be enlarged to
include the information on these variables. 
They can be represented by additional nodes in the graph with the restrictions 
that the subgraph on the non--dynamic nodes must be a DAG (or chain graph, cf.\
Gottard, 2002) and no directed
edges are allowed to point from processes to time--fixed covariates. 
A process being locally independent of a time--fixed covariate 
means that the intensity does not depend on this particular covariate given
all the other covariates and information on the past of all processes.
$\delta$--separation can still be applied to find further local independencies.

The nodes in a local independence graph do not necessarily have to stand for
only one mark (or the associated counting process), marks can be combined into 
one node as has been done in Examples 2 and 3 with `Health status' and 
`Tumour size'.
This might be of interest when there are logical dependencies. For example if
a particular illness is considered then the events
`falling ill' and `recovering' from this illness
are trivially locally dependent.
If one node is used to represent a collection of marks corresponding to a multivariate
subprocess  of the whole multivariate counting process then
an arrow into this node will mean that the intensity of at least one 
(but not necessarily all) of these marks depends on the origin of the arrow. 
However, some interesting information could be lost.
For instance, if events such as `giving birth to 1st child', `giving birth to 2nd child' etc.\
are considered, it might be relevant whether a woman is married or not
when considering the event of `giving birth to 1st child' but it might not be
relevant anymore when considering `giving birth to 2nd child.'

In many data situations the mark space is not finite, e.g.\ when measuring the 
magnitude of electrical impulse,  the amount of income in a new job or the dosage 
of a drug.
One could then discretise the mark space e.g.\ in `finding a well paid job'
and `finding a badly paid job'. However, it must be suspected that
too many of these types of events will generate too many logical dependencies
that are of no interest and will make the graphs crammed.

As we have seen in some of the examples it is sometimes sensible to consider
stopped processes in order to avoid having to represent
logical and uninteresting dependencies.  More generally one might want to relax 
in Definition \ref{def_li_mul} the requirement `for all $t\in {\cal T}$'
and instead consider suitably defined intervals based on stopping times.
For example, it might be the case that the independence structure is very different 
between the time of finishing education and starting the first job than before 
or after that. This deserves further investigation.

\section{Discussion and conclusions}\label{goodfor}

The main point of graphical models is that they allow certain
algebraic manipulations to be replaced by graphical ones. 
In the case of local independence graphs,  
we can read properties of intensity processes with respect to different,
in particular reduced filtrations, off the graph without the need to derive explicit
formulae for these intensities, and similarly we can read off relations among
subprocesses such as properties (\ref{total_fact}) and (\ref{part_fact}). 
This facilitates reasoning about complex dependencies, especially
in the face of unobservable information, 
and simplifies calculations by reducing dimensionality.

Clearly, it is tempting to interpret local independence graphs causally.
However, we regard causal inference as a topic of its own and it is not 
the aim of this paper to  go
into much detail in this respect, but for the following few comments.
Local independence graphs represent (in)dependencies in
$E(N_k(dt)|\calf_{t^-})$, where conditioning is on 
\em having observed \em $\calf_{t^-}$ and, as we saw in Example 2, 
it makes a difference
to what dependencies there are whether we condition on $\calf_{t^-}$
or different subsets (or even extensions) thereof.
Causal inference is about predicting $N_k(dt)$ after \em intervening \em in 
$\calf_{t^-}$, e.g.\ by modulating the times of the home visits to be once a week
in the home visits example. 
It is well known that conditioning on observation is not the same as
conditioning on intervention (``seeing'' and ``doing'' in  Pearl (2000)). 
Hence, without further assumptions, 
the arrows in local independence graphs do not necessarily represent
causal dependencies --- the intensity of an event being dependent on whether another
event has been observed before does not imply causation in the same way as 
correlation does not imply causation. 
Such further assumptions could be that
all `relevant' events (or processes) have been taken into account,
like originally proposed by Granger in order to justify the use of the term 
`causality' for what is now known as Granger-causality. 
E.g.\ if, in Figure \ref{dementia3}(b), we are satisfied that by including `Health'
all relevant processes have been taken into account, then we could say that homevisits
are indirectly causal for `Death'. Obviously, in this particular example,
there are many other relevant processes, like the occurrence of illnesses 
or  death of the partner, that might be relevant. 
However, 
the literature on (non-dynamic) graphical models and causality has
shown that causal inference is possible under weaker assumptions.
Analogous results based on local independence graphs would
require more prerequisites than we have given in this paper. 
Hence this is a topic for further research.
For non--graphical approaches to causal reasoning in a continuous time 
event history setting confer Eerola (1994), Lok (2001), Arjas and 
Parner (2004) and Fosen et al.\ (2004).

Another issue is the question of statistical inference for local 
independence graphs. This can be subdivided into (i) 
inference when
a graphical structure is given, e.g.\ from background knowledge, but we still 
want to quantify the strength of the dependencies, and (ii) finding
the graph from data if nothing about the local 
independence structure is known beforehand, which can be regarded as a particular
kind of model selection or search task. 
The former has partly been addressed in section \ref{sec_lik}, where more 
specific results will depend on the actual modelling assumptions about the
intensity processes which in turn will depend on the particular application. 
Estimation and testing within the class of Markov processes is
tackled in Didelez (2007).
More generally, local independence graphs can be combined with non--, semi-- or 
parametric methods but more research is required to investigate
how the graphical representation of the 
local independence structure can simplify inference in particular settings.
As to  model search, Nodelman et 
al.\ (2003) provide a first attempt, restricted to Markov processes, 
at exploiting the graphical structure  to find the
graph itself when it is not postulated based on
background knowledge. Clearly, generalisations would be desirable.

\begin{appendix}
\section{Appendix}

The appendix is targeted at proving Theorem \ref{theo_equiv}, but to do so we first give some
more results on the properties of local independence and $\delta$--separation 
which will be used in that proof. These 
are explored  along the lines of the graphoid axioms 
(Dawid, 1979; Pearl and Paz, 1987; Pearl, 1988; Dawid, 1998) 
which have been generalised to the asymmetric case by Didelez (2006).

\subsection{Properties of local independence}

\begin{prop} Properties of local independence \label{app_prop_li}\\
\em The following properties hold for local independence:
\begin{itemize}
\item[(i)] \em left redundancy: \em for all $A,B\subset V$: $A\li B|A$,
\item[(ii)] \em left decomposition: \em for all $A,B,C\subset V$ and $D\subset A$: if $A\li B|C$ then $D\li B|C$, 
\item[(iii)] \em left weak union: \em for all $A,B,C\subset V$ and $D\subset A$:
if $A\li B|C$ then $A\li B|(C\cup D)$  and\\
\em right weak union: \em for all $A,B,C\subset V$ and $D\subset B$:
if $A\li B|C$ then $A\li B|(C\cup D)$,
\item[(iv)] \em left contraction: \em \\
\hspace*{\fill} for all $A,B,C,D\subset V$: if $A\li B|C$ and $D\li B|(A\cup C)$ then $(A\cup D)\li B|C$,
\item[(v)] \em right intersection: \em \\
\hspace*{\fill} for all $A,B,C\subset V$: if $A\li B|C$ and $A\li C|B$ then $A\li (B\cup C)|(B\cap C)$.
\end{itemize}
\end{prop}
{\bf Proof:}\\ (i) Left redundancy holds since obviously the ${\cal F}_t^{A\cup
B}$--intensities of $\bfN_B$ are ${\cal F}_t^{A\cup B}$--measurable, i.e.\ if the past of
$\bfN_A$ is known, then the past of $\bfN_A$ is of course irrelevant.
\\
(ii) Left decomposition holds since the ${\cal F}_t^{A\cup B\cup C}$--intensities
$\lambda_k(t)$, $k\in B$, are ${\cal F}_t^{B\cup C}$--measurable by assumption so that the
same must hold for the ${\cal F}_t^{B\cup C\cup D}$--intensities $\lambda_k(t)$, $k\in B$,
for $D\subset A$.
\\
(iii) Left and  right weak union also trivially hold since adding information on the past of
components that are already uninformative (left) or included (right) does not change the 
intensity.\\
(iv) Left  contraction holds since we have that the ${\cal F}_t^{A\cup B\cup C\cup
D}$--intensities $\lambda_k$, $k\in B$, are by assumption ${\cal F}_t^{A\cup B \cup
C}$--measurable and these are again by assumption ${\cal F}_t^{B\cup C}$--measurable.\\
(v) 
The property of right intersection can be checked by noting that in the definition of local
independence the filtration w.r.t.\ which the intensity process should be
measurable is always generated at least by the process itself.  \bo\\

Note that left redundancy, left decomposition and left contraction imply that
\begin{eqnarray}
A\li B \mid C \quad \Leftrightarrow \quad A\backslash C \li B \mid C\label{disj1}.
\end{eqnarray}
It is also always true that $A\li B\mid C$ $\Rightarrow$ $A\li B\bsl C| C$, 
but we do not have equivalence here. \\

The following property will be important for the
equivalence of pairwise, local and global dynamic Markov  properties  (just like in the well--known 
case of  undirected conditional
independence graphs, cf.\ Lauritzen, 1996). 

\begin{prop} Left intersection for local independence\\
\em
Under the assumption of (\ref{filter})
local independence satisfies the following property called \em left intersection: \em for all
$A,B,C\subset V$
\bea
\mbox{ if }A\li B\mid C \mbox{ and } C\li B\mid A\mbox{ then }(A\cup C)\li B\mid (A\cap C)\nn   .
\eea
 \em
\end{prop}
{\bf Proof:}\\
Left intersection assumes that the ${\cal F}_t^{A\cup B\cup C}$--intensities $\lambda_k(t)$,
$k\in B$, are ${\cal F}_t^{B\cup C}$-- as well as ${\cal F}_t^{A\cup B}$--measurable.
With (\ref{filter}) we get that they are ${\cal F}_t^{B\cup (A\cap C)}$--measurable which yields
the desired result. \\

The following can be regarded as an alternative version of the above property of
left intersection. With (\ref{disj1}), left decomposition and left intersection 
we have that for disjoint $A,B,C,D\subset V$ 
\begin{eqnarray}
A\li B \mid (C \cup D) \mbox{ and } C\li B \mid (A \cup D) \Rightarrow
(A\cup C) \li B \mid D \label{inter1}
\end{eqnarray}
This follows from Corollary 4.3 of Didelez (2006).\\

The last property that I want to consider, the `right' counterpart of left decomposition given above, makes a
statement about the irrelevance of a process $\bfN_A$ after discarding part of the possibly \em
relevant \em information
$\bfN_{B\bsl D}$. If the irrelevance of $\bfN_A$ is due to knowing the past of
$\bfN_{B\bsl D}$ then it will not
necessarily be irrelevant anymore if the latter is discarded. 

\begin{prop} \label{assrdec} Conditions for right decomposition of local independence\\
\em
Consider a marked point process and assume that the cumulative counting process
$\sum N_k$ is non--explosive and that intensites exist.
Let $A,B,C\subset V$, $D\subset B$, with
$(B\cap A)\bsl (C\cup D)=\emptyset$.
The following property, called \em right decomposition, \em 
\[ A\li B\mid C \; \Rightarrow\; A\li D\mid C
\]
holds under the conditions that
\bea
B\li A\bsl(C\cup D)\mid(C\cup D) \nonumber
\eea
and
\bea
A\li \{k\}\mid C \cup B\, \mbox{ or } \, B\li \{k\}\mid (C \cup D\cup A) \nonumber
\eea
for all $k\in C\bsl D$.
\end{prop}

{\bf Proof:}\\
In this proof we proceed somewhat informally for the sake of simplicity. 
The formal proof is based on the results of Arjas et al.\ (1992)
and given in Didelez (2000, pp.\ 72).

Redefine $A^*=A\bsl (C\cup D)$, $B^*=B\bsl D$ and $C^*=C\bsl D$.
Then $A^*\cap B^*=\emptyset$ 
and with  (\ref{part_fact}) and (\ref{total_fact}), 
it can be shown that the assumptions of the present proposition imply that
\bea
N_D(t)\ind {\cal F}_{t^-}^{A^*} \mid {\cal F}_{t^-}^{B^*\cup C^*\cup D} \label{ind_1}
\eea
as well as
\bea
{\cal F}_t^{A^*}\ind {\cal F}_t^{B^*} \mid {\cal F}_t^{C\cup D}. \label{ind_2}
\eea
We want to show that the ${\cal F}_t^{A\cup C\cup D}$--intensity
$\tilde \lam_D(t)$ of $\bfN_D(t)$ is ${\cal F}_t^{C\cup D}$--measurable.
With the above and interpretation (\ref{compens}) we have
\bea
\tilde \lam_D(t)dt 
&=& E(N_D(dt)\mid \calf_{t^-}^{A\cup C\cup D}) =
E(N_D(dt)\mid \calf_{t^-}^{A^*\cup C^*\cup D})\nn \\
&=& E ( E(N_D(dt)\mid \calf_{t^-}^{A^*\cup B^*\cup C^*\cup D})\mid
\calf_{t^-}^{A^*\cup C^*\cup D}) \nn \\
&=& E ( E(N_D(dt)\mid \calf_{t^-}^{B^*\cup C^*\cup D})\mid
\calf_{t^-}^{A^*\cup C^*\cup D}) \nn \quad \quad\mbox{ using (\ref{ind_1})}\\
&=& E(N_D(dt)\mid \calf_{t^-}^{C^*\cup D})\nn \quad\quad\quad\quad\quad \mbox{ using (\ref{ind_2})}\\
&=& E(N_D(dt)\mid \calf_{t^-}^{C\cup D}),\nn 
\eea
as desired.

\subsection{Properties of $\delta$--separation} \label{app_delta}

For a general investigation of the properties of $\delta$--separation we need to complete 
Definition \ref{defi_delta_sep} by the case  that $A, B$ and $C$ are not disjoint: we then define that $C$
$\delta$--separates $A$ from $B$ if $C\bsl B$ $\delta$--separates $A\bsl (B\cup C)$ from $B$.
We further define that the empty set is always $\delta$--separated from $B$. Additionally, we define that the empty set
$\delta$--separates $A$ from $B$ if $A$ and $B$ are unconnected in $(G^B_{\aan(A\cup
B)})^m$.\\

It can be shown (Didelez, 2006) that $\delta$--separation satisfies the same properties 
as local independence given above in Proposition \ref{app_prop_li}
if we replace $A\li B|C$ by ``$C$ $\delta$--separates $A$ 
from $B$'' which we will write as $A\dsep B|C$. 
 In particular it satisfies left redundancy, left decomposition, let and right weak union, left and right
contraction as well as left and right intersection without requiring further asumptions. The property of
right decomposition holds for $\delta$--separation under conditions analogous to those 
in Proposition \ref{assrdec}. In particular we have the following special case of right decomposition
\bea
A\dsep B\mid C,\, D\subset B \Rightarrow A\dsep D\mid (C\cup B)\bsl D
\label{sprdec}
\eea
which is  Lemma 4.11 in Didelez (2006).\\

In addition, we want to show how $\delta$--separation can be read off a local independence graph in 
a different but equivalent way to Definition \ref{defi_delta_sep}.
We mention this firstly, because it will be more familiar to readers who use $d$--separation for DAGs 
(Pearl, 1988; Verma and Pearl, 1990), and secondly, because
some parts of the proof of Theorem \ref{theo_equiv} are easier to show using this alternative way of checking 
$\delta$--separation.

First, the different notions of paths and trails need to be made more stringent. 
Consider a directed or undirected graph $G=(V,E)$. An ordered $(n+1)$--tuple $(j_0, \ld, j_n)$
of distinct vertices is called an \em undirected path \em from
$j_0$ to $j_n$ if $\{j_{i-1},j_i\}\in E$ and a \em directed path \em if
$(j_{i-1},j_i)\in E$ for all $i=1,\ld, n$.
A (directed) path of length $n$ with $j_0=j_n$ is called a \em (directed) cycle. \em A subgraph
$\pi=(V',E')$ of $G$ with $V'=\{j_0,\ld, j_n\}$ and $E'=\{e_1, \ld, e_n\}\subset E$ is called a
\em trail \em between $j_0$ and $j_n$ if $e_i=(j_{i-1},j_i)$ or
$e_i=(j_{i},j_{i-1})$ or $e_i=\{j_{i},j_{i-1}\}$ for all $i=1,\ld, n$.
Further, for a directed graph 
we say that a trail between $j$ and $k$ is \em blocked by $C$ \em if it contains 
a vertex $\gamma$ 
such that either (i) directed edges of the trail do not meet head--to--head at
$\gamma$ and $\gamma\in C$, or (ii) directed edges of the trail meet 
head--to--head at $\gamma$ and $\gamma$ as well as all its descendants are no elements of $C$. 
Otherwise the trail is called \em active.  \em

\begin{prop} Trail condition for $\delta$--separation \label{protrail}\\
\em
Let $G=(V,E)$ be a directed graph and $A,B,C$ pairwise disjoint subsets of $V$. 
Define that any
\em allowed trail from $A$ to $B$ \em  contains no edge of the form
$(b,k), b\in B, k\in V\bsl B$.
For disjoint subsets $A,B,C$ of $V$, we have that $C$ $\delta$--separates
$A$ from $B$ if and only if all allowed trails from $A$ to $B$ are blocked by $C$.
\end{prop}
The proof is given in Didelez (2000, pp.22; cf.\ also Didelez, 2006).

\subsection{Proof of Theorem \ref{theo_equiv}}\label{app_proof}

It is easily checked that (\ref{glob}) $\Rightarrow$ (\ref{locmp}) $\Rightarrow$ (\ref{pair}):
First, pa$(k)$ always $\delta$--separates  $V\bsl (\mbox{pa}(k)\cup \{k\})$ from $\{k\}$ in
$G$, hence (\ref{locmp}) is just a special case of (\ref{glob}). 
Also, it is easy to see that
the equivalence of the pairwise and local dynamic Markov properties immediately
follows from left intersection assuming (\ref{filter}), left weak union and left decomposition.
Thus, the following proof considers situations where $A, B, C$ do not form a 
partition of $V$ or pa$(B)\subset\!\!\!\!\!/\; C$.
The structure of the proof corresponds to the one given by Lauritzen (1996, p.\ 34) for the
equivalence of the Markov properties in undirected conditional independence graphs. Due to the
asymmetry of local independence, however, this version
 is more involved. \\
Assume that (\ref{pair}) holds and that $C$ $\delta$--separates $A$ from $B$ in the local
independence graph. We have to show that $A\li B\mid C$, i.e.\ the ${\cal F}_t^{A\cup B\cup
C}$--intenities $\lambda_k(t)$, $k\in B$, are
$\calf_t^{B\cup C}$--measurable.
The proof is via backward induction on the number $|C| $ of vertices in the separating set. If
$| C| =| V| -2$ then both, $A$ and $B$ consist of only one element and (\ref{glob}) trivially
holds. If $| C| <| V| -2$ then either $A$ or $B$ consist of more than one element.
\\
Let us first consider the case that $A,B,C$ is a partition of $V$ and none of them is empty. If
$|A| >1$ let
$\alpha\in A$. Then, by left weak union and left decomposition 
of $\delta$--separation we have that
$C\cup (A\bsl\{\alpha\})$  $\delta$--separates  $\{\alpha\}$ from $B$, i.e.
\[
\{\alpha\}\dsep B\mid  C\cup (A\bsl\{\alpha\})
\]
and $C\cup \{\alpha\}$ $\delta$--separates $A\bsl \{\alpha\}$ from $B$ in $G$, i.e.
\[ A\bsl \{\alpha\}\dsep B\mid (C\cup \{\alpha\}).
\]
Therefore, we have by the induction hypothesis that
\[ \{\alpha\}\li B\mid C\cup (A\bsl\{\alpha\}) \mbox{ and } A\bsl \{\alpha\}\li B\mid (C\cup \{\alpha\}).
\]
From this it follows with the modified version of left intersection 
as given in (\ref{inter1})
(which can be applied because of the assumption that (\ref{filter}) holds) that $A\li B\mid C$ as
desired.\\ 
If $| B| >1$  we can show by a similar reasoning, applying  (\ref{sprdec})
to $\{\beta\}\in B$, that $A\li B\mid C$.\\
Let us now consider the case that $A,B,C\subset V$ are disjoint but no partition of $V$. First,
we assume that they are a partition of An$(A\cup B\cup C)$, i.e.\ that $A\cup B\cup C$ is an
ancestral set. Let $\gamma\in V\bsl(A\cup B\cup C)$, i.e.\ $\gamma$ is not an 
ancestor of $A\cup
B\cup C$. Thus, every allowed trail (cf.\ Proposition \ref{protrail}) from $\gamma$ to $B$
is blocked by
$A\cup C$ since any such trail includes an edge $(k,b)$ for
some $b\in B$ where no edges meet head--to--head in $k$ and $k\in A\cup
C$. Therefore, we get
\[
\{\gamma\}\dsep B\mid (A\cup C).
\]
Application of left contraction, weak union, and decomposition for
$\delta$--separation yields
\[
A\dsep B\mid (C\cup \{\gamma\}).
\]
It follows with the induction hypothesis that
\[
A\li B\mid (C\cup \{\gamma\}) \mbox{ as well as }\{\gamma\}\li B\mid (A\cup C).
\]
With left intersection as given by (\ref{inter1}) 
and left decomposition for local independence we get the desired result.\\
Finally, let $A,B,C$ be disjoint subsets of $V$ and $A\cup B\cup
C$  not necessarily an ancestral set. Choose $\gamma\in$ an$(A\cup B\cup
C)$ and define $\tilde G^B=G^B_{\aan(A\cup B\cup
C)}$. Since $A\gsep B\mid C$ in $(\tilde G^B)^m$ we know
from the properties of ordinary graph separation that
\bea
&(i)&\mbox{either }\{\gamma\}\gsep B\mid (A\cup C)  \mbox{ in }(\tilde G^B)^m \nn \\
&(ii)&\mbox{or } A\gsep \{\gamma\}\mid (B\cup C)  \mbox{ in }(\tilde G^B)^m.
\nn
\eea
(i) In the first case $\{\gamma\}\dsep B\mid (A\cup C)$ and it follows from left contraction that
\[
(A\cup \{\gamma\}) \dsep B\mid C.
\]
Application of left weak union and left decomposition yields $A\dsep B\mid (C\cup \{\gamma\})$.
With the induction hypothesis we therefore get
\bea
A\li B\mid (C\cup \{\gamma\}) \mbox{ and } \{\gamma\}\li B\mid (A\cup C).
\nn
\eea
Left intersection according to (\ref{inter1}) and left decomposition for local independence
yields $A\li B\mid C$.\\ 
(ii) The second case is the most complicated and the proof makes now use of
right decomposition for local independence under the conditions given in Proposition
\ref{assrdec}. First, we have from (ii) that $A\ind_g\{\gamma\}\mid B\cup C$ in
$(G_{\aan (A\cup B\cup C)})^m$ since the additional edges starting in $B$ can only yield
additional paths between $A$ and $\gamma$ that must be  intersected by $B$.
Since deleting further edges out of $\gamma$ does not create new paths, it holds
\[
A\dsep \{\gamma\}\mid (B\cup C).
\]
With $A\dsep B\mid C$, application of right contraction for $\delta$--separation yields $A\dsep (B
\cup \{\gamma\})\mid C$. Now, we can apply property
(\ref{sprdec}) to get $A\dsep B\mid  (C\cup \{\gamma\})$ from where it
follows with the induction hypothesis that
\[
A\li B\mid (C \cup \{\gamma\}) \mbox{ and } A\li \{\gamma\}\mid (B\cup C).
\]
With right intersection for local independence we get $A$ $\li$ $(B\cup
\{\gamma\})\mid C$.
In addition, $\{\gamma\}\li A\mid (B\cup C)$ by the same arguments as given above for $A\li
\{\gamma\}\mid  (B\cup C)$.
In order to apply Proposition \ref{assrdec} we still have to show that for all $k\in C$ either
$A\dsep \{k\}\mid (C \cup B\cup \{\gamma\})$ or $\{\gamma\}\dsep \{k\}\mid  (C\cup B\cup A)$
which by the induction hypothesis implies the corresponding local independencies. To see this,
assume that there exists a vertex $k\in C$ for which neither holds. With the trail condition
we then have that in $G_{\aan(A\cup B\cup C)}$ there exists an allowed trail from $A$ and
$\gamma$, respectively, to $k$ such that every vertex where edges do not meet head--to--head
are not in $(C\cup B\cup \{\gamma\})\bsl \{k\}$ and
$(C\cup B\cup A)\bsl \{k\}$,
respectively, and every vertex where edges meet head--to--head or some of their descendants
are in $(C\cup B\cup \{\gamma\})\bsl \{k\}$ respective $(C\cup B\cup A)\bsl \{k\}$. This would
yield a path between $A$ and $\gamma$ which is not blocked by $C\cup B$ (note that $k$ is a
head--to--head node on this trail) in $G_{\aan(A\cup B\cup C)}$. This in turn contradicts the
separation of $A$ and $\gamma$ by $B\cup C$ in $G^B_{\aan(A\cup B\cup C)}$ because the edges
starting in $B$ cannot contribute to this trail. Consequently we can apply right decomposition
and get the desired result.

\subsection{Proof of Theorem \ref{theo_cis}}\label{app_cis}

The factorisation (\ref{likfact}) implies that the marginal
likelihood for the marked point process discarding events not in 
An$(A\cup B\cup C)$, i.e.\ $\{
(T_s,E_s) \mid  s=1,2,\ld; E_s \in \cale_{\aan(A\cup B\cup C)} \}$, is given by 
\bea
L(t\mid  H_{t}^{\aan(A \cup B \cup C)}) = \prod_{k \in \aan(A \cup B\cup C) } L_k(t\mid 
H_{t}^{\cl(k)}) \nn
\eea
as none of the intensities of $N_k, k\in$An$(A\cup B\cup C)$ depend on 
$V\backslash$An$(A\cup B\cup C)$.
Hence, the likelihood may be written as product over factors that only depend on
cl$(k)$, $k\in$ An$(A\cup B\cup C)$. Let ${\cal C}=\{$cl$(k)\mid  k\in$
An$(A\cup B\cup C)\}$ be the set containing all such sets and
let $g_c(t|\cdot)$, $c\in {\cal C}$, be these factors. Then, we have
\bea
L(t\mid  H_{t}^{\aan(A\cup B\cup C)}) = \prod_{c\in {\cal C}} g_c(t\mid  H_{t}^{c}).  \nn
\eea
Further, by rearranging the factors
the sets in $\cal C$ can be taken to be the `cliques', 
i.e.\ the maximal fully
connected sets of nodes,  of the graph $(G_{\aan(A\cup B\cup C)})^m$. 
The above thus corresponds to the factorisation property of 
undirected graphs which in turn implies the
global Markov property for undirected graphs
(Lauritzen, 1996, p.\ 35). This means that when we have
a separation, like $C$ separating $A$ and $B$, in this graph
$(G_{\aan(A\cup B\cup C)})^m$ 
the corresponding conditional independence (\ref{total_fact}) holds,
which completes the proof.\\

\end{appendix}

\vspace{8mm}
\noindent{\bf Acknowledgements:}\\
I wish to thank Iris Pigeot for her persistent encouragement and advice, and
Niels Keiding for valuable comments and suggestions. 
Financial  support by the German Research Foundation (SFB 386) is also gratefully acknowledged.
The referees' comments have greatly helped to improve the paper.

\end{document}